\documentclass[11pt]{article}
\newcommand*{\mbq}{{\mathbb{Q}}}
\newcommand*{\noin}{\noindent}
\newcommand*{\zp}{{\mathbb{Z}_{p}}}
\newcommand*{\kar}{\longrightarrow}
\newtheorem{theorem}{Theorem}[section]
\newtheorem{lemma}[theorem]{Lemma}
\newtheorem{prop}[theorem]{Proposition}

\usepackage{eucal}
\usepackage{amscd}
\usepackage{amssymb}
\usepackage{amsmath}
\usepackage{amsfonts}
\usepackage{amsopn}
\usepackage{latexsym}
\begin{document}
\begin{center}
{\bf{\large On units generated by Euler systems}}\\
\vspace{.3cm}
{{\bf Anupam Saikia}}\\
{\it Department of Mathematics, IIT Guwahati,}\\
{\it Guwahati 781039.}\\
{\it Email: a.saikia@iitg.ernet.in}
\end{center}
\vspace*{.3cm} {\it Abstract}: In the context of cyclotomic
fields, it is still unknown whether there exist Euler systems
other than the ones derived from cyclotomic units. Nevertheless,
we first give an exposition on how norm-compatible units are
generated by any Euler system, following work of Coates. Then we
prove that the units obtained from Euler systems and the
cyclotomic units generate the same $\mathbb{Z}_{p}$-module for any
odd prime $p$. The techniques adopted for the Iwasawa theoretic
proof in latter part of this article originated in Rubin's work on
main conjectures of Iwasawa theory.

\section{Introduction}

Euler systems were introduced by Thaine and Kolyvagin. Later,
Rubin used Euler system of cyclotomic units and elliptic units to
prove the main conjecture of Iwasawa theory in various set-ups. In
[Co], Coates gave a definition of Euler systems in the context of
elliptic curves. His definition is somewhat stronger and different
from that of Rubin and Kolyvagin, but it fits more closely with
earlier work of Coates and Wiles ([CW 1], [CW 2]). In the context
of cyclotomic fields, the definition analogous to [Co] is as
follows. Let $S$ be any finite set of rational primes, always
containing 2. Let
\[W_{S}=\{ \zeta \in \bar{Q} : \zeta^{n}=1\;\;\mbox{with}\;(n,S)=1\}.\]

\noindent {\sc Definition :} An Euler system is a map $\phi :
W_{S} \longrightarrow \bar{Q}^{\times}$ which satisfies the
following axioms:
\begin{itemize}
\item (E 1) $\phi(\eta^{\sigma}) =
\phi(\eta)^{\sigma}\;\forall\;\sigma \in Gal(\bar{\mbq}/\mbq)$,
and $\phi(\eta^{-1})=\phi(\eta)$. \item (E 2) If $p$ is any
rational prime not in $S$, we have
\[\prod_{\zeta \in\mu_{p}}\phi(\zeta \eta)=\phi(\eta^{p}) \;\forall \eta \in W_{S}.\]
\item
(E 3) Let $p$ be any rational prime not in $S$. Then, for all
$\eta \in W_{S}$ of order prime to $p$, and all $\zeta \in
\mu_{p}$, we have
\begin{equation}
\phi(\zeta\eta)\equiv
\phi(\eta)\;\mbox{mod}\;{\mathfrak{p}}\;\forall \mathfrak{p}|p.
\label{mod}
\end{equation}
\end{itemize}
\noin[Here, $\mathfrak{p}$ is a prime ideal over $p$ in any field
that contains $\phi(\zeta\eta)$].\\

\noindent The basic example of an Euler system is essentially the
classical system of cyclotomic units. We briefly describe them.
Let $\Omega$ denote the non-zero integers $a_{1},\ldots a_{r}$ and
integers $n_{1},\ldots , n_{r}$ such that $\sum_{j=1}^{r}
n_{j}=0$. Let $\lambda_{\Omega}(T)$ be given by
\[\lambda_{\Omega}(T) = \prod_{j=1}^{r}(T^{-a_{j}}-T^{a_{j}})^{n_{j}}.\]
Let $S$ be the set consisting of 2 and all primes $q$ such that
$q$ divides at least one of the $a_{i}$. Let $\phi_{\Omega} :
W_{S} \longrightarrow \bar{Q}^{\times}$ be defined as \\
$(i)\;
\phi_{\Omega}(\eta)=\lambda_{\Omega}(\eta)\;\mbox{for}\;\eta \in
W_{S}$
and $\eta \neq 1$.\\
$(ii)\;\phi_{\Omega}(1)=\lim\limits_{T\rightarrow
1}\lambda_{\Omega}(T) =
\prod\limits_{j=1}^{r} a_{j}^{n_{j}}$.\\

\noindent It is easy to check that $\phi_{\Omega}$ forms an Euler
system. Given an Euler system $\phi :  W_{S} \longrightarrow
\bar{\mathbb{Q}}^{\times}$, there are two ways of constructing new
ones:
\begin{itemize}
\item If $n$ is any non-zero integer, $\phi\circ n$ is again an
Euler system provided we enlarge $S$ to include primes dividing
$n$. \item If $\xi$ is a primitive $h$-th root of 1, define
$\phi_{\xi}(\eta) =\prod\limits_{\tau}\phi(\eta \xi^{\tau})$ where
$\tau$ runs over the elements of the Galois group
Gal($\mathbb{Q}(\xi)/\mathbb{Q}$). Here, we enlarge $S$ by
including all the primes that divide $h$.
\end{itemize}
It is still unknown whether there exist Euler systems attached to
cyclotomic fields other than those mentioned above.

\section{Euler systems generate global units}

In this section, we will explain how Euler systems attached to
cyclotomic fields generate norm-compatible global units in the
cyclotomic tower. The next three propositions are consequences of
the axioms (E 1), (E 2) and (E 3), and are cyclotomic analogues of
results in [Co]. For each $m \geq 1$, let $\mu_{m}$ denote the
group of $m$-th roots of unity. If $\eta \in \mu_{m}$, with
$(m,S)=1$, then (E 1) shows that $\phi(\eta)\in
\mbq(\mu_{m})^{+}\subset \mbq(\mu_{m})$, where $L^{+}$ denotes the
maximal real subfield of $L$. Let $\mbq_{m}$ denote the cyclotomic
field $\mbq(\mu_{m})$ and $H_{m}$ denote the maximal real subfield
of $\mbq_{m}$, i.e., $H_{m}=\mbq(\mu_{m})^{+}$.
\begin{prop}
Let $\eta$ be an element of $\mu_{m}$ with $(m,S)=1$. Let ${p}$ be
any prime with $({p},{m})=1$, and ${p}\not \in S$, Then, we have
\begin{equation}
\label{frob} N_{H_{mp}/H_{m}}\phi(\zeta
\eta)=\phi(\eta)^{\mbox{Frob}_{p}-1} \;\; \forall \;\zeta \in
\;\mu_{p},\; \eta \neq 1.
\end{equation}
Here, $Frob_{{p}}$ denotes the Frobenius element of
 ${p}$ in the Galois group of $H_{m}$ over $\mathbb{Q}$, which is
 unramified at $p$.
\end{prop}
Proof: By axiom E1 and E2, \begin{align*} \scriptstyle{
N_{H_{mp}/H_{m}}\phi(\zeta \eta)} &\scriptstyle{\; =\;
\prod\limits_{\sigma\in Gal(H_{mp}/H_{m})}\phi(\zeta
\eta)^{\sigma}\; =\; \prod\limits_{\sigma\in Gal(\mbq_{mp}/\mbq_{m})}\phi(\zeta^{\sigma}\eta^{\sigma})}\\
&\scriptstyle{ = \;\frac{\prod\limits_{\xi \in
\mu_{p}}\phi(\xi\eta)}{\phi(1.\eta)}\;=\;
\frac{\phi(\eta^{p})}{\phi(\eta)}\;=
\;\phi(\eta)^{Frob_{p}-1}.\quad \square}
\end{align*}
\begin{prop}
Let $\eta$ be any element of $\mu_{m}$ with $(m,S)=1$. Let $p$ be
any prime with $(p,m)=(p,S)=1$. For each $n\geq 0$, let
$\zeta_{n}$ be a primitive ${p}^{n+1}$-root of 1 such that
$\zeta_{n+1}^{p}=\zeta_{n}$. Then the sequence
$\phi(\zeta_{n}\eta^{\mbox{Frob}_{p}^{-n}})(n=0, 1, \ldots )$, is
norm compatible in the tower $H_{mp^{\infty}}$ over $H_{mp}$.
\end{prop}
Proof: By axioms E1 and E2,
\begin{align*}
\scriptstyle{ N_{H_{mp^{n+1}}/H_{mp^{n}}}\phi(\zeta_{n}
\eta^{Frob_{p}^{-n}})\;}\;&\scriptstyle{=\;\prod\limits_{\sigma\in
{Gal}(H_{mp^{n+1}}/H_{mp^{n}})}
\phi(\zeta_{n}\eta^{Frob_{p}^{-n}})^{\sigma} \;=\;
\prod\limits_{\sigma\in Gal(\mbq_{mp^{n+1}}/\mbq_{mp^{n}})}
\phi(\zeta_{n}^{\sigma}\eta^{{Frob_{p}^{-n}}})}\\
& \scriptstyle{\;= \;\prod\limits_{\xi \in
\mu_{p}}\phi(\xi\zeta_{n}\eta^{{Frob}_{p}^{-n}})\;=\;
\phi(\zeta_{n}^{p}\eta^{p{{Frob}_{p}^{-n}}})\;=\;
\phi(\zeta_{n-1}\eta^{{Frob}_{p}^{-(n-1)}}).\qquad \square}
\end{align*}
\begin{prop} For all $\eta \in W_{S}$ with $\eta \neq
1,\;\phi(\eta)$ is a unit.
\end{prop}
In order to prove the above proposition, we need the following
lemma:
\begin{lemma}
Let $p$ be any prime and $K/\mbq$ be a finite extension. Let
$\alpha \in K$ be a universal norm in the tower
$K(\mu_{p^{\infty}})$. Then every prime ideal in the factorization
of $\alpha$ divides $p$.
\end{lemma}
\noin Proof: Let $\mathfrak{q}$ be a prime not dividing $p$ that
occurs in the factorization of $\alpha$. Now, the Galois group
Gal($K(\mu_{p^{\infty}})/K$) is a subgroup of
$\mathbb{Z}_{p}^{\times}$. After finitely extending $K$, we can
assume that the Galois group is $\zp$. Since $\mathfrak{q}$ is
unramified in the tower, it is enough to show that its
decomposition group in non-trivial. Then, $\mathfrak{q}$ will be
an inert prime, and infinite power of $\mathfrak{q}$ will divide
the universal norm $\alpha$, which is absurd. If the decomposition
group is trivial, it will imply that the residue fields in the
tower $K(\mu_{p^{\infty}})/F$ are finite. But, if $\zeta$ and
$\xi$ are $p$-power roots of unity, then $\zeta \not \equiv \xi$
modulo $\mathfrak{q}$ unless $\mathfrak{q}$ lies
above $p$. \\

\noin Proof of proposition 2.3: Let $\eta$ be a primitive root of
order $m=p^{n+1}b$, where $(p,b)=1$. Let us denote the number
field $\mbq(\eta)^{+}$ by $K$. We know that $\phi(\eta)\in K$. Let
$\mathfrak{q}$ be a prime ideal in $K$ that divides $\phi(\eta)$.
By proposition 2.2, $\phi(\eta)$ is a universal norm in the tower
$\mbq(\mu_{p^{\infty}}\eta)^{+}$ over $K=\mbq(\eta)^{+}$. Clearly,
$\phi(\eta)^{2}\in K$ is a universal norm in the tower
$K(\mu_{p^{\infty}})$ over $K$. Hence, $\mathfrak{q}$ must divide
$p$ by the above lemma. If $b$ is not 1 we are through, as
$\mathfrak{q}$ must also divide any prime factor of $b$. \\

\noin Now consider the case $b=1$. Then $\eta=\zeta_{n}$, and by
proposition 2.2,
$N_{H_{p^{n+1}}/H_{p}}\phi(\zeta_{n})=\phi(\zeta_{0})$. Now,
\begin{equation*}
\begin{split}
N_{H_{p}/\mbq}\phi(\zeta_{0})& =
\big(N_{\mbq(\mu_{p})/\mbq}\phi(\zeta_{0})\big)^{\frac{1}{2}}\\
 & =
\big(\prod_{\sigma\in {Gal}(\mbq(\mu_{p})/\mbq)}\phi(
\zeta_{0})^{\sigma}\big)^{\frac{1}{2}}\\
&= \big(\prod_{\sigma\in {Gal}(\mbq(\mu_{p})/\mbq)}\phi(
\zeta_{0}^{\sigma})\big)^{\frac{1}{2}}\\
&  =
\Bigg(\frac{\prod\limits_{\xi \in
\mu_{p}}\phi(\xi)}{\phi(1)}\Bigg)^{\frac{1}{2}}\\
& =\Big(\frac{\phi(1)}{\phi(1)}\Big)^{\frac{1}{2}} =\pm 1.
\end{split}
\end{equation*}
Noting that $K=\phi(\zeta_{n})^{+}$ is totally ramified over
$\mbq$ at $p$, we have only one prime $\mathfrak{q}$ of $K$ above
$p$. But now $(\phi(\eta))=\mathfrak{q}^{r}$ for some integer $r$,
and $r$ has to be zero as the norm of $\phi(\eta)$ is $\pm 1$.
Thus, $\phi(\eta)$ is a global unit in the ring of integers of
$K$. \;\;$\quad\square$\\

\noin Thus, we can conclude that Euler systems attached to
cyclotomic fields generate norm-compatible global units.

\section{Statement of main result} Let $p$ be an odd prime, and
$\zeta_{n}$ be a fixed $p^{n+1}$-th root of unity such that
$\zeta_{n+1}^{p} = \zeta_{n} \forall n\geq 0$. Let $F_{n}$ denote
the number field $\mathbb{Q}(\zeta_{n})^{+}$. Let
\[{\cal{H}}(p) =\;\; \mbox{the set of Euler systems}\;\;\phi :
W_{S} \longrightarrow \bar{\mathbb{Q}}^{\times} \;\; \mbox{such
that}\;\;p \not\in S.\] By proposition 2.3, above,
$\phi(\zeta_{n})$ is a global unit in $F_{n}$ for any $\phi\in
{\cal{H}}(p)$. By proposition 2.2, \begin{equation} \label{normm}
N_{F_{n+1}/F_{n}}\phi(\zeta_{n+1})=\phi({\zeta_{n}}).
\end{equation}
Let us define
\[{\cal{E}}_{n}=\{\phi(\zeta_{n})\mid \phi \in {\cal{H}}(p)\}.\]
Let $E_{n}$ and $C_{n}$ denote respectively the global units and
the cyclotomic units in $F_{n}$. \noin Let $\mathfrak{p}_{n}$ be
the unique maximal ideal of $F_{n}$ above p. We denote the
completion of $F_{n}$ at $\mathfrak{p}_{n}$ by $\Phi_{n}$. Let
$U_{n}$ be the principal local units of $\Phi_{n}$, i.e., the
local units in $\Phi_{n}$ congruent to 1 mod $\mathfrak{p}_{n}$.
Let ${E}_{n,1}$, ${C}_{n,1}$ and ${\mathcal{E}}_{n,1}$ denote
respectively the subgroup of principal units in $E_{n}$, $C_{n}$
and ${\cal E}_{n}$. These subgroups can be canonically embedded in
$U_{n}$. Let $\tilde{E}_{n,1}$, $\tilde{C}_{n,1}$ and
$\tilde{\cal{E}}_{n,1}$ be respectively the closure of
${E}_{n,1}$, ${C}_{n,1}$ and ${\mathcal{E}}_{n,1}$ in $U_{n}$. Let
$\bar{E}_{n}=\zp\otimes E_{n}$, $\bar{C}_{n}=\zp\otimes C_{n}$ and
$\bar{\cal{E}}_{n}=\zp\otimes {\cal{E}}_{n}$. By Leopoldt's
conjecture, which is proved to be true for the abelian extensions
$F_{n}$ of $\mbq$, we have $\bar{E}_{n}=\tilde{E}_{n,1}$,
$\bar{C}_{n}= \tilde{C}_{n,1}$ and
$\bar{\cal{E}}_{n}=\tilde{\cal{E}}_{n,1}$. Hence we have a natural
inclusion
\[\bar{C}_{n} \subset \bar{\cal E}_{n}\subset \bar{E}_{n}. \]
 The main result of
this paper is the theorem below:
\begin{theorem} The $\mathbb{Z}_{p}$-module generated by
the global units derived from Euler systems attached to $p$-power
cyclotomic fields is the same as the $\mathbb{Z}_{p}$-module
generated by the cyclotomic units. In other words,
${\bar{\cal{E}}}_{n}=\bar{C}_{n}$.
\end{theorem}

\noin Note that the index of $\bar{C}_{n}$ in $\bar{E}_{n}$ equals
$h_{p,n}$, where $h_{p,n}$ denotes the $p$-part of the class
number of $F_{n}=\mbq(\zeta_{n})^{+}$. Vandiver has conjectured
that $h_{p,0}=1$ (which is equivalent to saying that $h_{p,n}=1$
for all $n=0,1,2,\ldots$). If one can show that
$[\bar{E}_{0}:\bar{\cal{E}}_{0}]=1$, then it will imply Vandiver's
conjecture by virtue of theorem 3.1. Of course, even to attempt
this approach to Vandiver's conjecture, one would certainly
require an Euler system which is not derived from cyclotomic
units.\\

\noin Another interesting point to note here is the connection of
theorem 3.1 to Greenberg's conjecture. Greenberg's conjecture is
equivalent to the statement that the only universal norms in
${\bar{E}_{0}}$ is the group ${\bar{C}_{0}}$. Since the values of
Euler systems are universal norms, theorem 3.1 gives evidence for
Greenberg's conjecture. One can raise the question whether the
only universal norms in each ${\bar{E}_{n}}$ are those coming from
Euler systems.\\

\noin We will prove the above theorem by establishing a relation
involving Iwasawa modules, and then by descent.

\section{Iwasawa theoretic set-up}

\noindent Let us consider the infinite extension
\[F_{\infty} = \cup_{n\geq 0}F_{n},\]
where $F_{n}=\mathbb{Q}(\zeta_{n})^{+}$ and $\zeta_{n}$ is a
primitive $p^{n+1}$-th root of unity such that
$\zeta_{n+1}^{p}=\zeta_{n}$. Let us define the Galois groups
\[G_{n}= G(F_{n}/\mathbb{Q}),\qquad \Gamma_{n}= G(F_{n}/F_{0}).\]
Let $G_{\infty}$ be the Galois group of $F_{\infty}$ over
$\mathbb{Q}$ and $\Gamma$ be the Galois group
Gal($F_{\infty}/F_{0}$). Clearly,
\[G_{\infty}= G(F_{\infty}/\mathbb{Q}) = \Delta \times \Gamma ,\;\; \mbox{where}\;
\Delta \simeq G(F_{1}/\mathbb{Q}).\]\\
\noin We have the following
field diagram:
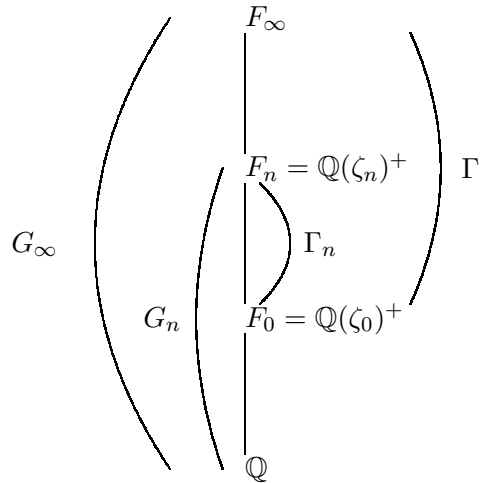
\begin{center}
\setlength{\unitlength}{1mm}
\begin{picture}(80,75)
\put(40,70){\makebox(0,0)[l]{$F_{\infty}$}}
\put(40,50){\makebox(0,0)[l]{$F_{n}=\mathbb{Q}(\zeta_{n})^{+}$}}
\put(40,30){\makebox(0,0)[l]{$F_{0}=\mathbb{Q}(\zeta_{0})^{+}$}}
\put(40,10){\makebox(0,0)[l]{$\mathbb{Q}$}}
\put(40,52){\line(0,1){16}} \put(40,32){\line(0,1){16}}
\put(40,12){\line(0,1){16}} \qbezier(30,70)(10,40)(30,10)
\put(12,40){\makebox(0,0){$G_{\infty}$}}
\qbezier(37,50)(30,30)(37,10) \put(29,30){\makebox(0,0){$G_{n}$}}
\qbezier(62,68)(70,50)(62,32) \put(70,50){\makebox(0,0){$\Gamma$}}
\qbezier(42,48)(50,40)(42,32)
\put(50,40){\makebox(0,0){$\Gamma_{n}$}}
\end{picture}
\end{center}
Let $R_{n}$ be the group ring of $G_{n}$ with coefficients in
$\mathbb{Z}_{p}$. These group rings form an inverse system under
the canonical maps from $R_{m}$ to $R_{n}$. We define
\[R_{\infty} = \lim_{\substack{\longleftarrow \\ n} }\mathbb{Z}_{p}[G_{n}].\]
A $R_{\infty}$-module $N$ is called a torsion $R_{\infty}$-module
if it is annihilated by a non-zero-divisor in $R_{\infty}$. If $N$
is a finitely generated torsion $R_{\infty}$-module, then there is
an injective $R_{\infty}$-module homomorphism
\begin{equation}
\label{st0} \bigoplus_{i=1}^{r}R_{\infty}/g_{i}R_{\infty}
\hookrightarrow N
\end{equation}
with finite cokernel. The elements $g_{i}$ are not uniquely
determined, but the ideal $\prod\limits_{i}g_{i}R_{\infty}$ is. We
call the ideal $\prod\limits_{i}g_{i}R_{\infty}$ the {\it
characteristic ideal} of $N$ and denote it by char{$(N)$}. The
characteristic ideal is multiplicative in exact sequence: if
$0\longrightarrow N' \longrightarrow N \longrightarrow N''
\longrightarrow 0$ is an exact sequence of torsion
$R_{\infty}$-modules then
\[\mbox{char}(N)=\mbox{char}(N')\mbox{char}(N'').\]
From now on, $M$ will always denote a fixed power of $p$. We
denote the group ring of $G_{n}$ with coefficients in
$\mathbb{Z}_{p}$ modulo $M$ by $R_{n,M}$, i.e.,
\[R_{n,M}=(\mathbb{Z}/M\mathbb{Z})\;[G_{n}] = R_{n}/MR_{n}.\]
We denote the $p$-part of the ideal class group of $F_{n}$ by
$A_{n}$. These groups form an inverse system under the norm maps
and we denote the inverse limit by $A_{\infty}$.\\

Now, the $\mathbb{Z}_{p}$-modules $U_{n}$, $\bar{C}_{n}$,
$\bar{\cal E}_{n}$ and $\bar{E}_{n}$ defined in the previous
section are equipped with an $R_{n}$-module structure. They form
an inverse system of $R_{n}$-modules with respect to the norm
maps. Thus, we define the inverse limits
\[U_{\infty} = \lim_{\substack {\longleftarrow \\ n}} U_{n},\;\;E_{\infty} =
\lim_{\substack {\longleftarrow \\ n}}
\bar{E_{n}},\;\;{\cal{E}}_{\infty}= \lim_{\substack
{\longleftarrow \\ n}}\bar{{\cal{E}}}_{n},\;\;C_{\infty}=
\lim_{\substack {\longleftarrow \\ n}}\bar{C_{n}}.\] These inverse
limits have the natural structure of a $R_{\infty}$-module. We
will first determine a relation between the Iwasawa modules
${\cal{E}}_{\infty}$ and $C_{\infty}$. Then, we will descend to
the $n$-th layer.
\begin{prop} The characteristic ideal of $A_{\infty}$ contains
the characteristic ideal of $E_{\infty}/{\mathcal{E}}_{\infty}$.
\end{prop}
\noindent {\it Proof of theorem 3.1 assuming proposition 4.1 : }
Our main result follows easily from proposition 4.1. We have the
following exact sequence of $R_{\infty}$-modules
\begin{align}
0 \longrightarrow {\cal{E}}_{\infty}/{C}_{\infty}\longrightarrow &
E_{\infty}/C_{\infty}\longrightarrow E_{\infty}/{\cal{E}}_{\infty}
\longrightarrow 0.\nonumber\\
\Rightarrow \;\mbox{char}( E_{\infty}/C_{\infty}) = & \mbox{char}(
E_{\infty}/{\cal E}_{\infty}) \mbox{char}({\cal
E}_{\infty}/C_{\infty}). \label{rr}
\end{align}
By ``main conjecture'' of Iwasawa theory for cyclotomic fields, we
have
\begin{equation}
 \mbox{char}(A_{\infty}) = \mbox{char}( E_{\infty}/C_{\infty})
\label{ss}
\end{equation}
Then proposition 4.1 combined with (\ref{rr}) and (\ref{ss}) imply
that $\mbox{char}({\cal E}_{\infty}/C_{\infty}) = R_{\infty}$.
That tells us that ${\cal E}_{\infty}/C_{\infty}$ is a finite
$R_{\infty}$ submodule of $ U_{\infty}/C_{\infty}$. However, it is
well-known that $ U_{\infty}/C_{\infty}\simeq R_{\infty}/(g)$. For
instance, the results in [Sa] gives an explicit proof of this
fact. Hence, $ U_{\infty}/C_{\infty}$ has no non-trivial finite
$R_{\infty}$-submodule. Thus, ${\cal E}_{\infty}/ C_{\infty}=0$,
and $(C_{\infty})_{\Gamma_{n}}=({\cal{E}}_{\infty})_{\Gamma_{n}}$.
However, by the $R_{\infty}$-module structure of $C_{\infty}$, we
know that $(C_{\infty})_{\Gamma_{n}}=\bar{C}_{n}$. Moreover,
$({\cal{E}}_{\infty})_{\Gamma_{n}}\twoheadrightarrow
\bar{\cal{E}}_{n}$. Thus, the canonical injection $\bar{C}_{n}\kar
\bar{\cal{E}}_{n}$ is also surjective, and we have our main
result.\\


\noin Our approach to proposition 4.1 will be as follows. We will
use the units generated by Euler systems to construct new elements
(called Kolyvagin class), which factorize non-trivially. We will
then determine the factorization of these new elements
(proposition 6.2). This gives us a systematic way of obtaining
relations in the ideal class group. These relations will be recast
in Iwasawa theoretic set-up, and we will obtain a suitable ideal
which annihilates the Iwasawa module of class groups. This
approach originated in Rubin's work. In the remaining sections, we
will closely follow Rubin's proof of the main conjectures of
Iwasawa theory.

\section{Kolyvagin class}

 From now on, we fix $n$ and refer
to $F_{n}$ simply as $F$, dropping the subscript. In this section,
we will construct elements in $F^{\times}$ using the global units
derived from any Euler system. The factorization of these elements
can be easily determined, which will be shown in the following
section. Let $M$ be a fixed power of $p$. Let ${\cal S}_{M} $ be
the set of square-free integers $s$ such that each prime factor
$q$ of $s$ splits in $F/\mathbb{Q}$ and $q\equiv 1
\;\mbox{mod}\;{M}$. For the rest of this paper, $q$ will always
denote a rational prime in ${\cal{S}}_{M}$, and $\mathfrak{q}$
will be a primes of $F$ above $q$. Let $\eta_{q}$ be a fixed
primitive $q$-th root of 1. We write $F(q)$ for the field
$F(\eta_{q})$ and $G(q)$ for the Galois group $G(F(\eta_{q})/F)$.
Clearly, $\mathfrak{q}$ is totally ramified in $F(q)/F$ and the
ramification index is $(q-1)$. Suppose $\sigma_{q}$ in $G(q) $
sends $\eta_{q}$ to $\eta_{q}^{t}$, where $t$ is a primitive root
mod $q$. Then $G(q)$ is cyclic and
generated by $\sigma_{q}$.\\

\noin As in [Ru 4], let us now define the following operators:
\[D_{q}= \sum_{i=1}^{q-2}i \sigma_{q}^{i},\;\;D_{s}=
\prod_{q \mid s}D_{q},\;\;\mbox{and}\;\;N_{q}=\sum_{i=0}^{q-2}
\sigma_{q}^{i}.\] It is easily seen that
\begin{equation}
\label{norm} (\sigma_{q}-1)D_{q} = (q-1-N_{q}).
\end{equation}
 From
(\ref{frob}) of pp. \pageref{frob}, it follows that
\begin{equation}
\label{frob1}
\begin{split}
N_{q}D_{r}\phi(\zeta_{n}\eta_{rq})\;& = \; D_{r} N_{F(rq)/F(r)}\phi(\zeta_{n}\eta_{rq})\\
 & =  \;(\mbox{Frob}_{q}-1)D_{r}\phi(\zeta_{n}\eta_{r}).
\end{split}
\end{equation}

\begin{prop} $D_{s}\phi(\zeta_{n}\eta_{s})$ is an element of
$\big((F(s)^{\times})/(F(s)^{\times})^{M}\big)^{G(s)}$. In other
words, $[(\sigma -1)D_{s}\phi(\zeta_{n}\eta_{s})]^{\frac{1}{M}}$
is a well defined element of $F(s)^{\times}$ for all $\sigma$ in
$G(s)$.
\end{prop}
{\sc Proof} : We use induction on the number of primes dividing
$s$. Suppose $q|s$ and $ s=qr$. Then
\begin{align*}
(\sigma_{q}-1)D_{s}\phi(\zeta_{n}\eta_{s}) & =
(\sigma_{q}-1)D_{q}D_{r}\phi(\zeta_{n}\eta_{s})\\
& =
(q-1)D_{r}\phi(\zeta_{n}\eta_{s})/(\mbox{Frob}_{q}-1)D_{r}\phi(\zeta_{n}\eta_{r})
\;\; {\text{(by (\ref{norm}) and (\ref{frob1}))}}.
\end{align*}
Since $q$ is in ${\cal{S}}_{M},\; M|(q-1)$. As $\mbox{Frob}_{q}
\in G(r)$, the induction hypothesis implies that the denominator
in the last expression above is in $(F(r)^{\times})^{M}$.
Therefore,
\[(\sigma_{q}-1)D_{s}\phi(\zeta_{n}\eta_{s}) \in (F(s)^{\times})^{M}.\]
Since the $\sigma_{q}$ generate $G(s)$, this completes the proof
of the proposition.
\hspace*{1cm}$\square$ \\

\noindent As $\sigma$ runs over the elements of $G(s)$, $\sigma
\mapsto [(\sigma -1)D_{s}\phi(\zeta_{n}\eta_{s})]^{\frac{1}{M}}$
gives an element of $H^{1}\big(G(s), F(s)^{\times}\big)$. By
Hilbert 90, this cohomology group is trivial. Therefore, there is
an element $\beta_{s,\phi}$ in $F(s)^{\times}$ such that
\begin{equation}
\label{koly} [(\sigma
-1)D_{s}\phi(\zeta_{n}\eta_{s})]^{\frac{1}{M}}=(\sigma-1)\beta_{s,\phi}.
\end{equation}
Clearly, $\beta_{s,\phi}$ is unique up to multiplication by an
element of
$F^{\times}$. We can now make the following definition.\\

\noindent {\sc Definition :} {\sl For} $s \in {\cal S}_{M} $, {\sl
Kolyvagin class is defined as}
\begin{equation}
\label{kol1} \kappa_{\phi,M}(s) =
\frac{D_{s}\phi(\zeta_{n}\eta_{s})}{\beta_{s,\phi}^{M}} \in
F^{\times}/( F^{\times})^{M},
\end{equation}
{\sl where} $\phi$ {\sl is an Euler system in} ${\cal{H}}(p)$ {\sl
and}
$\beta_{s,\phi}$ {\sl is given by (\ref{koly})}.\\

\section{Factorization of Kolyvagin class}

\noin In the previous section we constructed certain elements,
called Kolyvagin classes, in the $F^{\times}$ modulo $M$th powers,
where $F=\mbq(\zeta_{n})^{+}$. Here we will describe how one can
determine the factorization of those elements. Such a
factorization should be seen as a relation in the ideal class
group of $F$. This process will lead to construction of an
annihilator of
the class group of $F$ in the group ring $\frac{\mathbb{Z}}{M\mathbb{Z}}[Gal(F/\mbq]$. \\

 \noindent Let ${\cal O}_{F}$ be the ring of
integers of $F$ and
\[ I_{F} = I = \bigoplus \mathbb{Z}\mathfrak{q} \]
be the group of fractional ideals of $F$ written additively. Let
\[ I_{F,q} = I_{q} =  \bigoplus_{{\mathfrak{q}|q}} \mathbb{Z}\mathfrak{q}. \]
 For any $x \in F^{\times}$, let $(x)\in I $ be the principal ideal generated by
$x$, and $(x)_{q},\;[x]_{M}$, and $[x]_{q,M}$ the projections of
$(x)$ to $I_{q},\;I/MI$, and $I_{q}/MI_{q}$ respectively. When
there is no ambiguity, we drop the subscript $M$ and simply write
$[x]$ or $[x]_{q}$. Note that $[x]$ and $[x]_{q}$ are well defined
for $x\in F^{\times}/ (F^{\times})^{M}$. The next two propositions
are cyclotomic analogues of lemma 13
and theorem 14 in [Co]: \\

\noindent
\begin{prop} There is a Galois equivariant isomorphism
\[\lambda_{q}\;: \big( {\mathcal{O}}_{F}/q {\mathcal{O}_{F} }\big)^{\times}/
\Big( \big( {\mathcal{O}}_{F} /q
{\mathcal{O}}_{F}\big)^{\times}\Big)^{M} \longrightarrow
I_{q}/MI_{q}. \]
\end{prop}
{\sc Proof} : Let $\tilde{\mathfrak{q}}$ be the unique prime of
$F(q)$ above the prime $\mathfrak{q}$ of $F$ and
$\pi(\tilde{\mathfrak{q}})$ be a local parameter at
$\tilde{\mathfrak{q}}$. The residue fields of $q,\;\mathfrak{q}$
and $\tilde{\mathfrak{q}}$ will be denoted by $k(q), \;
k(\mathfrak{q})$ and $k(\tilde{\mathfrak{q}})$ respectively. As
$q$ splits in $F$ and $\mathfrak{q}$ is totally ramified in
$F(q)/F$, the residue fields are all isomorphic. We have an
isomorphism
\[G(q) \longrightarrow k(\tilde{\mathfrak{q}})^{\times},\qquad \sigma \mapsto
\pi({\tilde{\mathfrak{q}}})^{1-\sigma}
\;\mbox{mod}\;\tilde{\mathfrak{q}}.\] Note that since $G(q)$ is
the inertia group in $F(q)/F$, the above isomorphism does not
depend on the choice of the parameter $
\pi({\tilde{\mathfrak{q}}})$. If $ \sigma_{q}$ maps to
$\gamma(\tilde{\mathfrak{q}})$ under the above isomorphism, then
clearly  $\gamma(\tilde{\mathfrak{q}})$ is a generator of
$k(\tilde{\mathfrak{q}})^{\times}$. By our identification above,
$\gamma(\tilde{\mathfrak{q}})$ can be regards as a generator of
$k(\mathfrak{q})^{\times}$. For any $w \in \big
({\mathcal{O}}_{F}/q {\mathcal{O}}_{F}\big)^{\times}$, we have
\[ w \equiv \gamma(\tilde{\mathfrak{q}})^{a({\mathfrak{q}})}\;\mbox{mod}\;
\mathfrak{q}\;\;\mbox{for some
integer}\;a(\mathfrak{q})\;\;\mbox{mod}\;(q-1).\] Let us define
\[\lambda_{q}(w) = \sum_{\mathfrak{q} \mid q}(a(\mathfrak{q})\; \mbox{mod}\;M)\;
\mathfrak{q}.\] Galois equivariance and surjectivity follow
easily. Since both sides have the same cardinality, we have a
Galois equivariant isomorphism.
\hspace*{1cm}$\square$ \\
\noindent {\sl Note} : Since $q$ splits in $F$, we have a Galois
equivariant map
\[\bar{\lambda}_{q} : \big( {\mathcal{O}}_{F}/q {\mathcal{O}}_{F}\big)^{\times}/
\Big( \big( {\mathcal{O}}_{F}/q {\mathcal{O}}_{F}\big)^{\times}
\Big)^{M} \longrightarrow \mathbb{Z}/M\mathbb{Z}\; [G]\] given by
\begin{equation}
\label{barl} \bar{\lambda}_{q}(w)\mathfrak{q} = \lambda_{q}(w),
\end{equation}
where we fix a $\mathfrak{q}$ above $q$.\\

\noindent
\begin{prop}
 For $qs \in {\mathcal{S}}_{M}$, we have
\begin{align*}
{}&(i)\; \big[\kappa_{\phi,M}(s)\big]_{q} = 0.\\
{}&(ii)\; \big[\kappa_{\phi,M}(sq)\big]_{q} =
\lambda_{q}(\kappa_{\phi,M}(s)).
\end{align*}
\end{prop}

\noindent {\sc Proof} : Recall that
\[\kappa_{\phi,M}(s) = \frac{D_{s}\phi(\zeta_{n}\eta_{s})}{\beta_{s,\phi}^{M}}
\in F^{\times}/( F^{\times})^{M}.\] Since
$\phi(\zeta_{n}\eta_{s})$ is a unit in $F(s)$, so is
$D_{s}\phi(\zeta_{n}\eta_{s})$. Hence, the ideal generated by
$\kappa_{\phi,M}(s)$ in $F(s)$ is determined by
${\beta_{s,\phi}}$.
Note that $q$ does not divide $s$, because $qs$ is a square-free integer. \\

\noindent (i) $q$ is unramified in $F(s)$, and
$\beta_{s,\phi}^{M}$
is an $M$-th power in $F(s)$. \\

\noindent (ii) Let $\tilde{\cal{Q}}$ be a prime of $F(sq)$ above
the prime $\mathfrak{q}$ of $F$. The ramification index of
$\tilde{{\cal{Q}}}$ in $F(sq)/F$ is $(q-1)$. By definition of
$\kappa_{\phi,M}(sq)$,
\begin{equation}
\label{val} v_{\mathfrak{q}}(\kappa_{\phi,M}(sq)) = -
\frac{M}{q-1}v_{\tilde{{\cal{Q}}}}(\beta_{sq,\phi}).
\end{equation}
Since $F(sq)/F(q)$ is unramified at ${\tilde{\cal{{Q}}}}$, the
local parameter $ \pi(\tilde{\mathfrak{q}})$ at the prime ideal
$\tilde{\mathfrak{q}}$ of $F(q)$ is also a local parameter at
$\tilde{\cal{{Q}}}$. We have the following diagram of fields and
prime ideals:
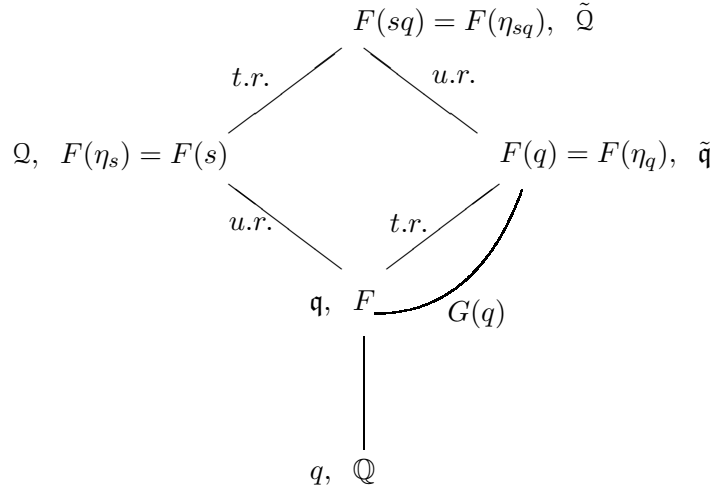
\begin{center}
\setlength{\unitlength}{1.5mm}
\begin{picture}(80,53)
\put(41,10){\makebox(0,0)[r]{$q,\;\;\mathbb{Q}$}}
\put(41,25){\makebox(0,0)[r]{$\mathfrak{q},\;\;F$}}
\put(52,38){\makebox(0,0)[l]{$F(q)=F(\eta_{q}),\;\;\tilde{\mathfrak{q}}$}}
\put(39,50){\makebox(0,0)[l]{$F(sq)=F(\eta_{sq}),\;\;\tilde{\mathcal{Q}}$}}
\put(28,38){\makebox(0,0)[r]{${\mathcal{Q}},\;\;F(\eta_{s})=F(s)$}}
\qbezier(41,24)(50,24)(54,35) \put(50,24){\makebox(0,0){$G(q)$}}
\put(44,32){\makebox(0,0){$t.r.$}}
\put(30,32){\makebox(0,0){$u.r.$}}
\put(48,45){\makebox(0,0){$u.r.$}}
\put(30,45){\makebox(0,0){$t.r.$}} \put(40,12){\line(0,1){10}}
\put(42,28){\line(4,3){10}} \put(38,28){\line(-4,3){10}}
\put(50,40){\line(-4,3){10}} \put(28,40){\line(4,3){10}}
\end{picture}
\end{center}
Let $v_{\tilde{{\cal{Q}}}}(\beta_{sq,\phi}) = b$. Then,
\begin{align*}
\beta_{sq,\phi}\;&=
\;\pi(\tilde{\mathfrak{q}})^{b}u,\;\;\text{where
$u \in F(sq)$ is prime to $\tilde{\cal{Q}}$.} \\
\Rightarrow \;\beta_{sq,\phi}^{1-\sigma_{q}} \; & = \;
\pi(\tilde{\mathfrak{q}})^{({1-\sigma_{q}}){b}}u^{1-\sigma_{q}} \\
                               \; & \equiv \; \gamma(\tilde{\mathfrak{q}})^{b}\;\;
                               \text{mod $\tilde{\cal{{Q}}}$ (as $\sigma_{q}$
                               acts trivially on $u$ modulo $\tilde{\cal{{Q}}}$).}
\end{align*}
From (\ref{koly}),
\begin{eqnarray}
                     \beta_{sq,\phi}^{\sigma_{q} -1}
&     =  & \big[( {\sigma_{q}
-1})D_{sq}\phi(\zeta_{n}\eta_{sq})\big]^{
           \frac{1}{M}}\nonumber\\
&     =  & \big[
D_{s}\phi(\zeta_{n}\eta_{sq})^{(q-1)}/D_{s}\phi(\zeta_{n}
           \eta_{s})^{\mbox{Frob}_{{q}} - 1} \big]^{\frac{1}{M}}
            \;\;[\mbox{by (\ref{norm}) and (\ref{frob1})}] \nonumber\\
&     =  &  \big[
D_{s}\phi(\zeta_{n}\eta_{sq})^{(q-1)}/(\beta_{s,\phi}^M)^{
           {\mbox{Frob}_{{q}} - 1}} \big]^{\frac{1}{M}}\;\;\mbox{[by (\ref{koly})]}
           \nonumber\\
&     \equiv  &  \big[ D_{s}\phi(\zeta_{n}\eta_{s})^{(q-1)}/
(\beta_{s,\phi}^M)^{{\mbox{Frob}_{{q}} - 1}}
\big]^{\frac{1}{M}}\;\;\mbox{mod}\;\;
            \tilde{\cal{Q}}\;\;\mbox{[by (\ref{mod})]}  \nonumber\\
& \equiv &  \big[ D_{s}\phi(\zeta_{n}\eta_{s})/
\beta_{s,\phi}^{M}\big]^{\frac{q-1}{M}}
\;\;\mbox{mod}\;\tilde{\cal{Q}}\;\;
[\mbox{since Frob}_{q}(x) \equiv x^{q}\;\mbox{mod}\;\tilde{\cal{Q}}]\nonumber\\
            \label{kol3}\Rightarrow \gamma(\tilde{\mathfrak{q}})^{-b}
& \equiv &  \kappa_{\phi,M}(s)^{\frac{q-1}{M}} \;\;\mbox{mod}\;\;
            \tilde{\cal{Q}}\;\;[\mbox{by definition of}\;\;\kappa_{\phi,M}(s)].
\end{eqnarray}
If $\lambda_{q}( \kappa_{\phi,M}(s)) = \sum_{\mathfrak{q} \mid
q}(a(\mathfrak{q})\; \mbox{mod}\; M)\mathfrak{q}$,
\begin{align*}
\kappa_{\phi,M}(s) \; & \equiv \; \gamma
(\tilde{\mathfrak{q}})^{a(\mathfrak{q})}\;\;\mbox{mod}\;{\mathfrak{q}}\;\;
\text{[by definition of $\lambda_{q}$]}\\
\Rightarrow \; {\gamma(\tilde{\mathfrak{q}})}^{-b}\; & \equiv
\;{\gamma(\tilde{\mathfrak{q}})}^{a(\mathfrak{q})\frac{q-1}{M}}\;\;
\mbox{mod}\;{\mathfrak{q}} \;\;\text{[by (\ref{kol3})]}\\
\intertext{Since $\gamma(\tilde{\mathfrak{q}})$ is a generator of
$k(\mathfrak{q})^{\times}$, we have} \frac{q-1}{M}a(\mathfrak{q})
\;
& \equiv \; -b \;\mbox{mod}\;{(q-1)}.\\
\Rightarrow \;v_{\mathfrak{q}}(\kappa_{\phi,M}(sq))\;&\equiv
\;a(\mathfrak{q})\;\;\text{mod}\;M \;\;\text{by (\ref{val})}
\qquad\square
\end{align*}

\section{Kolyvagin sequence}

\noindent By standard Iwasawa theory, $(A_{\infty})_{\Gamma_{n}}=
A_{n}$ (see [Wa] for a proof). Since $A_{n}$ is finite,
$A_{\infty}$ is a finitely generated torsion $R_{\infty}$-module
by Nakayama's lemma. By structure theory of $R_{\infty}$-modules,
there is an exact sequence of $ R_{\infty}$-modules
\begin{equation}
\label{st} 0 \longrightarrow
\bigoplus_{i=1}^{r}R_{\infty}/(f_{i})\longrightarrow
A_{\infty}\longrightarrow D \longrightarrow 0,
\end{equation}
where $D$ is finite. By definition,
$\mbox{char}(A_{\infty})=\prod\limits_{i}f_{i}R_{\infty}$. Let
$y_{i}\in A_{\infty}$ be the image of $ 1 \in R_{\infty}/(f_{i})$
and $A_{\infty}^{0} =  \sum\limits_{i=1}^{r}R_{\infty}y_{i}$.
Suppose $J$ is the annihilator of $D=A_{\infty}/A_{\infty}^{0}$.
Then $J$ is of finite index in $R_{\infty}$. Let $M$ be a fixed
power of $p$.\\

\noin From (\ref{st}), we obtain an exact sequence
\begin{equation}
D^{\Gamma_{n}} \longrightarrow \bigoplus_{i=0}^{r}R_{n}/(f_{i})
\longrightarrow A_{n}\longrightarrow D_{\Gamma_{n}}
\longrightarrow 0 . \label{Nn}
\end{equation}
It follow that $R_{n}/\mbox{char}(A_{\infty})R_{n}$ is finite. Let
\[ N_{n} = \mid
A_{n}\mid.\mid\frac{R_{n}}{\mbox{char}(A_{\infty})R_{n}}\mid,
\qquad M' = MN_{n}.\]

\noindent {\sc Definition :} {\sl Let} $0 \leq k \leq r $. {\sl A
Kolyvagin sequence of length} $k$ {\sl is a} $k$-{\sl tuple}
$\mathfrak{Q}=(\mathfrak{q}_{1},\mathfrak{q}_{2},\ldots
,\mathfrak{q}_{k})$ {\sl of primes in} $F$ {\sl such that}
\begin{itemize}
\item
{\sl the} $\mathfrak{q}_{i}$ {\sl lie above distinct rational
primes in} $ {\cal S}_{M} $, {\sl and}
 \item
$\mbox{ Frob}_{\mathfrak{q}_{i}} = \sigma_{i}|_{L_{n}}$, {\sl
where} $\sigma_{i}\in G(L_{\infty}/F_{\infty})$ {\sl corresponds
to} $y_{i}\in A_{\infty}$ {\sl under Artin reciprocity} [\it where
$y_{i}$ are as defined immediately after (\ref{st})].
\end{itemize}

\noindent {\sl For a Kolyvagin sequence} $\mathfrak{Q}$, {\sl we
define the square-free integer} $s(\mathfrak{Q})$ {\sl as}
\[s(\mathfrak{Q}) = \prod_{i=1}^{k}q_{i}, \qquad \mbox{where}\;\;q_{i}
=\mathfrak{q}_{i}|_{\mathbb{Q}}.\]

\section{Key proposition}

\noin In this section, we work out the details of the proof of
proposition 4.1. It is a simplified version of Rubin's arguments
in (see [Ru 1], [Ru 2], [Ru 3], and [Ru 4]) together with ideas of
Coates [Co].\\

\noin Let $ \varPi(k,n,M)$ be the set of all Kolyvagin sequence of
length $k$. Let $\Psi(k,n,M)$ be the ideal in $R_{n,M} =
R_{n}/MR_{n}$ generated by
\[ \{ \psi(\kappa_{\phi,M}(s(\mathfrak{Q})))\mid \mathfrak{Q} \in \varPi(k,n,M),\psi
\in \mbox{Hom}_{R_{n}} ( R_{n,M} \kappa_{\phi,M}(s(\mathfrak{Q})),
R_{n,M})\}.\] In order to prove 4.1, we need the following key
proposition.

\begin{prop} $J \Psi(k,n,MN_{n})R_{n,M} \subset f_{k+1}\Psi(k+1,n,M)$.
\end{prop}
The following result is from Rubin [Ru 4]:
\begin{lemma}
 Let $B$ be a $p$-torsion free finitely generated $\mathbb{Z}_{p}[G]$-module where
 $G$ is a finite abelian group. If $f\in \mathbb{Z}_{p}[G]$ is not a zero divisor,
 $b\in B $, and
\[\{ \psi(b) : \psi \in \mbox{Hom}_{\mathbb{Z}_{p}[G]}(B,\mathbb{Z}_{p}[G])\}
\subset f \mathbb{Z}_{p}[G],\] then $b \in fB$.
\end{lemma}
\noindent Proposition 4.1 can be deduced from proposition 8.1 and lemma 8.2 as follows : \\

\noindent For $k=0$, $\varPi(0,n,M)$ has just the empty sequence
and
\[ \{ \psi(\kappa_{\phi,M}(1))\; \mbox{mod}\;M\;|\;\;\psi \in
\mbox{Hom}_{R_{n}}( \bar{E}_{n}, R_{n})\} \subset \Psi(0,n,M).\]
In other words,
\[\{ \psi(\phi(\zeta_{n}))\; \mbox{mod}\;M\;|\;\;\psi \in
\mbox{Hom}_{R_{n}}( \bar{E}_{n}, R_{n})\}\subset \Psi(0,n,M).\] If
we use proposition 3.4 recursively, we obtain
\begin{align*}
{} & J^{r} \Psi(0,n,MN_{n}^{r})R_{n,M} \subset
\mbox{char}(A_{\infty})\Psi(r,n,M)\subset
\mbox{char}(A_{\infty})R_{n,M}\\
\Rightarrow\; & J^{r} \psi(\phi(\zeta_{n}))R_{n,M} \subset
\mbox{char}(A_{\infty})R_{n,M}
\quad \forall \psi \in \mbox{Hom}_{R_{n}}( \bar{E}_{n}, R_{n}) \quad \forall M \\
\Rightarrow \;& J^{r} \psi(\phi(\zeta_{n}))R_{n} \subset
\mbox{char}(A_{\infty})R_{n}.
\end{align*}
Substituting $\bar{E}_{n}$ for $B$ in the lemma above, we get
\[ J^{r} (\phi(\zeta_{n}))\subset  \mbox{char}(A_{\infty})\bar{E}_{n}.\]
Taking inverse limit as $n$ goes to infinity, we find that
\[ J^{r} {\cal E}_{\infty}\subset   \mbox{char}(A_{\infty})E_{\infty}.\]
There is an obvious exact sequence
\[ 0 \longrightarrow {\mathcal{E}}_{\infty}/ J^{r} {\mathcal{E}}_{\infty}\longrightarrow
E_{\infty}/ J^{r} {\mathcal{E}}_{\infty}\longrightarrow
E_{\infty}/{\mathcal{E}}_{\infty} \longrightarrow 0. \] Since
$J^{r}$ is an ideal of finite index in $R_{\infty}$,
${\mathcal{E}}_{\infty}/ J^{r} {\mathcal{E}}_{\infty}$ is finite
and so its characteristic ideal is trivial. Thus,
\begin{equation}
\mbox{char}( E_{\infty}/{\mathcal{E}}_{\infty})= \mbox{char}(
E_{\infty}/J^{r} {\cal E}_{\infty}). \label{tt}
\end{equation}
The exact sequence of $R_{\infty}$-modules
\[ 0 \longrightarrow \mbox{char}(A_{\infty}) E_{\infty}/ J^{r} {\cal E}_{\infty}
\longrightarrow  E_{\infty}/ J^{r}
{\cal{E}}_{\infty}\longrightarrow E_{\infty}/
\mbox{char}(A_{\infty}){ E}_{\infty}\longrightarrow 0 \] implies
that
\begin{equation}
\mbox{char}( E_{\infty}/\mbox{char}(A_{\infty}){ E}_{\infty}) |
\mbox{char}( E_{\infty}/ J^{r} {\cal E}_{\infty}). \label{uu}
\end{equation}
It is clear from the structure theorem of $R_{\infty}$-modules
that $\mbox{char} (A_{\infty}) \mid
\mbox{char}(E_{\infty}/\mbox{char}(A_{\infty}){ E}_{\infty})$.
Therefore, (\ref{tt}) and (\ref{uu}) imply that
\[\mbox{char}(A_{\infty}) | \mbox{char}( E_{\infty}/{\cal E}_{\infty}).\]

\noindent In order to facilitate the proof of proposition 8.1,
we need the following lemma.\\
\begin{lemma}
Let $0 \leq k \leq r $ and
\[\mathfrak{Q} = ( \mathfrak{q}_{1}, \mathfrak{q}_{2},\ldots \mathfrak{q}_{k+1})\in
\varPi(k+1,n,MN_{n})\]
 be a  Kolyvagin sequence. Let $q$ be the rational prime below $\mathfrak{q}_{k+1} =
  \mathfrak{q} $ and $q_{i}$ be the rational prime below $\mathfrak{q}_{i}$ for
  $1\leq i \leq k$. Let $s = \prod\limits_{i=1}^{k}q_{i}$. Then there is a Galois
   equivariant map
\[\tilde{\psi} : R_{n,M} \kappa_{\phi,M}(sq) \longrightarrow   R_{n,M}, \]
 such that for any $\rho \in J$, we have
\[ \rho \bar{\lambda}_{q}(\kappa_{\phi,M'}(s)) \equiv f_{k+1}\tilde{\psi}
(\kappa_{\phi,M}(sq))\;\;\mbox{mod}\;M.\]
\end{lemma}

\noindent{\sc Proof} : Let $\mathcal{C}_{i}$ be the class of $
\mathfrak{q}_{i}$ in $A_{n} = A$. Let $\bar{A}_{n}$ be the
quotient of $A_{n}$ by the $R_{n}$-submodule generated by classes
of $ {\mathcal{C}}_{1}, {\mathcal{C}}_{2},\ldots
{\mathcal{C}}_{k}$. Let $A_{\infty}^{k} =
\sum_{i=0}^{k}R_{\infty}y_{i} \subset A_{\infty}^{0}$. In the
exact sequence (\ref{st}) of pp. \pageref{st}, the ideal $J$ of
$R_{\infty}$ annihilates $D$. Let $J_{n}$ be the image of $J$ in
$R_{n}$. From the exact sequence (\ref{Nn}) of pp. \pageref{Nn},
it is clear that $J_{n}$ annihilates the kernel of the map
\[  (A_{\infty}^{0}/A_{\infty}^{k}) \otimes R_{n} \longrightarrow \bar{A}_{n}.\]
Let $\mathcal{A }$ be the annihilator of the image
$\bar{{\mathcal{C}}_{k}}$ of the ideal class of $\mathfrak{q}$ in
$\bar{A}_{n}$. The annihilator of the class of $ \mathfrak{q}$ in
$ (A_{\infty}^{0}/A_{k}^{0}) \otimes R_{n} $ is $f_{k+1}R_{n}$. So
we have
\begin{equation}
\label{jay}
 J_{n}{\mathcal{A}} \subset {f_{k+1}}R_{n}.
\end{equation}
By proposition 6.2 (i), we have
\begin{equation}
( \kappa_{\phi,M'}(sq)) =
[\kappa_{\phi,M'}(sq)]_{q}\;\;\mbox{mod}\;(M'I,\;
R_{n}\mathcal{C}_{1}+ R_{n} \mathcal{C}_{2}+\ldots + R_{n}
\mathcal{C}_{k} ). \label{vv}
\end{equation}
Since $|A|$ divides $M'=MN_{n}$, $[ {\kappa}_{\phi,M'} (sq)]_{q}
\in {\mathcal{A}} R_{n,M'}{\mathfrak{q}}$, and by (\ref{jay}), $
\rho [{\kappa}_{\phi,M'}(sq) ]_{q} \in f_{k+1} R_{n,M'}$. Since
$f_{k+1}$ is not a zero divisor , we have a well defined map
\begin{align*}
 f_{k+1}^{-1} : f_{k+1} R_{n,M'} &\longrightarrow R_{n,M}, \\
                f_{k+1}g  & \mapsto \; h, \;\;\;\;\text{where $h$ is given by}\\
      f_{k+1} g  & = \; f_{k+1} h + MN_{n} \tilde{h}\\
  (\Rightarrow g \; =  \; h + M.f_{k+1}^{-1} N_{n} \tilde{h}, \;\;
  \mbox{note that}\;\;f_{k+1}\;\;\mbox{divides}\;\;N_{n}.)
\end{align*}
We define a map
\begin{align}
\tilde{\psi} : R_{n,M} \kappa_{\phi,M}(sq) \;& \longrightarrow \;
R_{n,M},\;\;\text{by} \nonumber\\
\tilde{\psi}(\kappa_{\phi,M}(sq))\mathfrak{q} \;& =\;f_{k+1}^{-1}
\rho [\kappa_{\phi,M'}(sq)]_{q} \label{psi}
\end{align}
and extend by linearity to the whole of
$R_{n,M}\kappa_{\phi,M}(sq)$. We have to verify that
$\tilde{\psi}$ is well-defined. Let $ \kappa_{\phi,M}(sq)^{\tau} =
\alpha^{M}$ for some $\tau \in R_{n}$. We want to show that
\[\tau\tilde{\psi}(\kappa_{\phi,M}(sq)) \in M R_{n,M'}.\]
By (\ref{vv}),
\[ (\alpha)  = [\alpha ]_{q}\;\;\mbox{mod}\;( N_{n}I,\;
R_{n}\mathcal{C}_{1}+ R_{n} \mathcal{C}_{2}+\ldots + R_{n} \mathcal{C}_{k} ).\]
Recall that $M'=MN_{n}$, where $N_{n}$ is divisible by the order
of $A_{n}$. Therefore, $[\alpha ]_{q} \in {\mathcal{A}}
R_{n,M'}{\mathfrak{q}}$ and by (\ref{jay}),
\begin{align*}
\rho[\alpha ]_{q} \in f_{k+1}R_{n,M'}\mathfrak{q}\\
\tilde{\psi} ( \kappa_{\phi,M}(sq)^{\tau})\mathfrak{q}
{}& = \tau f_{k+1}^{-1}\rho [\kappa_{\phi,M}(sq)]_{q}\\
{}& = f_{k+1}^{-1} \rho M [\alpha]_{q}  \subset M R_{n,M'}.
\end{align*}
Hence $ \tilde{\psi}$ is well defined. Now, proposition 2.3 (ii),
(\ref{psi}) and (\ref{barl}) imply that
\[f_{k+1}\tilde{\psi}(\kappa_{\phi,M}(sq))\equiv \rho \bar{\lambda}_{q}
(\kappa_{\phi,M'}(s))\;\;\mbox{mod}\;M.\hspace*{1cm}\square \]\\

\noindent {\sc Proof of Proposition} 8.1 : Let us consider any
Kolyvagin sequence of length $k$ in $\varPi(k,n,MN_{n})$, say
$\mathfrak{Q}=(\mathfrak{q}_{1}, \mathfrak{q}_{2},\ldots
\mathfrak{q}_{k})$ and let $s = s(\mathfrak{Q})$. Suppose $\psi$
is an arbitrary element in $ \mbox{Hom}_{R_{n}}(
R_{n,M'}\kappa_{\phi,M'}(s), R_{n,M'})$. We want to show that
\[\rho \psi(\kappa_{\phi,M'}(s)) R_{n,M} \subset f_{k+1}
\Psi (k+1,n,M) \quad
\forall \rho \in J. \]
We will extend the above Kolyvagin sequence to one of length $k+1$.
We apply Cebotarev density theorem in the following way: \\

\noindent Let
\[ W = R_{n,M'}{ \kappa}_{\phi,M'} (s)
\subset F^{\times}/{( F^{\times})}^{M'} \hookrightarrow
F(\mu_{M'})^{\times}/{( F(\mu_{M'})^{\times})}^{M'}. \] See the
proof of proposition 15.47 in [Wa] for the second injection above.
Let $F' = F(\mu_{M'})$, $L= L_{n}$ (the maximal unramified abelian
$p$-extension of $F$) and $H =  F(\mu_{M'}, W^{\frac{1}{M'}})$. We
have the following diagram of fields:
 \begin{center}
\setlength{\unitlength}{1mm}
\begin{picture}(80,50)
\put(48,45){\makebox(0,0)[l]{$H=F(\mu_{M'},W^{\frac{1}{M'}})$}}
\put(48,23){\makebox(0,0)[lb]{$F'=F(\mu_{M'})$}}
\put(48,3){\makebox(0,0)[lb]{$F$}} \put(48,9){\line(-5,3){20}}
\put(50,9){\line(0,1){12}} \put(50,29){\line(0,1){12}}
\put(23,24){\makebox(0,0){$L$}}
\end{picture}
\end{center}
We have a Kummer pairing
\[ G(H/F') \times W \longrightarrow \mu_{M'},\]
and $G(F'/\mathbb{Q})$-isomorphism
\[ G(H/F')\stackrel{\sim}{\longrightarrow}\mbox{Hom}(W,\mu_{M'}).\]
Complex conjugation acts trivially on the ideals of the real field
$F$, therefore it acts trivially on $G(L/F)$ by global class field
theory. It clearly implies that complex conjugation acts trivially
on $G(LF'/F')$. However, complex conjugation acts non-trivially on
Hom($W, \mu_{M'}$), and hence on $G(H/F')$. Thus, $LF' \cap H =
F'$. Then ramification consideration tells us that $L\cap H = F$.
If we ignore the Galois structure, then we have
\begin{equation}
\label{gal} G(H/F') \simeq \mbox{Hom}(W,\mathbb{Z}/M'\mathbb{Z}).
\end{equation}
Let
\[ \tau : (\mathbb{Z}/M'\mathbb{Z})[G] \longrightarrow
\mathbb{Z}/M'\mathbb{Z},\quad \sum a_{g}g \longmapsto a_{1}.\]
We can define a map $\tau \circ \psi$ by composition:
\[\tau \circ \psi  : W \stackrel{\psi}{\longrightarrow} R_{n,M'}
= (\mathbb{Z}/M'\mathbb{Z})\;[G] \stackrel{\tau}{\longrightarrow}
\mathbb{Z}/M'\mathbb{Z}. \]
Let $\gamma$ in $G(H/F') \subset G(H/F)$ correspond to $\tau \circ
\psi$ under (\ref{gal}). Let $\sigma$ in $G(LH/F)$ be such that
\[ \sigma \mid_{L} = \sigma_{k+1} \mid_{L},
\qquad \sigma \mid_{H} = \gamma. \]
It is possible to find such a $ \sigma $ because $L\cap H = F$. By
Cebotarev Density Theorem, there are infinitely many  degree 1
primes in $F$ such that the corresponding Frobenius in $G(LH/F)$
is in the same conjugacy class as $\sigma$. We pick one such
$\mathfrak{q}$ which is unramified in $H/F$. Since only finitely
many primes are ramified in the finite extension $H/F$, such a
choice is possible. Let $q$ be the rational prime below
$\mathfrak{q}$. As $\mathfrak{q}$ is an unramified prime of degree
1, $q$ splits in $F$. Hence, $N_{F/\mathbb{Q}}\mathfrak{q} = q$.
Since $\sigma$ fixes $F(\mu_{M'})$, $\mbox{Frob}_{\mathfrak{q}}$
fixes the residue field of  $F(\mu_{M'})$ at $\mathfrak{q}$.
Therefore,
\[N_{F/\mathbb{Q}}\mathfrak{q} = q \equiv 1 \;\mbox{mod}\;{M'}.\]
Thus, $q$ is a prime in ${\cal{S}}_{M}$ (recall that $M' =
MN_{n}$). We can now extend $\mathfrak{Q}$ to a Kolyvagin sequence
$\mathfrak{Q}'=(\mathfrak{q}_{1}, \ldots \mathfrak{q}_{k},
\mathfrak{q})$ of length $k+1$. Clearly, $s(\mathfrak{Q}')=sq$. By
lemma 3.6, there is a Galois equivariant map
\[\tilde{\psi} : R_{n,M}\kappa_{\phi,M}(sq) \longrightarrow R_{n,M},\]
such that
\begin{equation}
\label{final} \rho \bar{\lambda}_{q}(\kappa_{\phi,M'}(s)) \equiv
f_{k+1}\tilde{\psi}(\kappa_{\phi,M}(sq))\;\;\mbox{mod}\;M.
\end{equation}

\noindent Let $\bar{\mathfrak{q}}$ be a prime of $H$ above
$\mathfrak{q}$ such that $\mbox{Frob}_{\bar{\mathfrak{q}}}=\sigma
$. Let $w\in W $. Then $ v_{\bar{\mathfrak{q}}}(w) \equiv
\;0\;\mbox{mod}\;M'$. Now,
\begin{eqnarray*}
\tau\circ\psi(w) = 0
& \Leftrightarrow & \gamma(w^{\frac{1}{M'}}) = w^{\frac{1}{M'}}\nonumber\\
& \Leftrightarrow &
\mbox{Frob}_{\bar{\mathfrak{q}}}(w^{\frac{1}{M'}})
                    = w^{\frac{1}{M'}}\\
& \Leftrightarrow & w\;\;\mbox{is a}\;\;M'-{\mbox{th power
mod}}\;\bar{\mathfrak{q}}\cap F= \mathfrak{q}.\nonumber
\end{eqnarray*}

\noindent Suppose $\bar{\lambda}_{q}(w) = \sum a_{g} g $. By
definition of $\bar{\lambda}_{q}$, $ \tau\circ
\bar{\lambda}_{q}(w)=0 $ iff $w$ is  $M'$-th power mod
$\mathfrak{q}$. Thus,
\[\tau\circ\bar{\lambda_{q}} (w)=0 \Leftrightarrow w\;\;\mbox{is a}
\;\;M'{\mbox{-th power mod}}\;\;\mathfrak{q}\Leftrightarrow \tau\circ\psi(w) = 0. \]
Then,
\[ \psi = u \bar{\lambda}_{q}, \qquad u\in(\mathbb{Z}/M'\mathbb{Z})^{\times}.\]
(The above statement can be easily proved, as shown in lemma 15.49 of [Wa]).\\
From (\ref{final}), it is now obvious that
\[ \rho \psi(\kappa_{\phi,M'}(s)) \equiv u \rho
\bar{\lambda}_{q}(\kappa_{\phi,M'}(s))  \equiv u f_{k+1}
\tilde{\psi}(\kappa_{\phi,M}(sq)) \;\mbox{mod}\;M.\]
It is now clear that
\[ J \Psi(k,n,MN_{n}) \subset f_{k+1} \Psi(k+1,n,M).\qquad\square
\]\\

\noin{\it Acknowledgement:} I would like to thank Prof. John
Coates and the anonymous referee for their valuable suggestions.

\end{document}